\documentclass[12pt]{article}

\begin{document}

\title{Statistically dual distributions 
and estimation of the parameters}

\author {S.I. Bityukov$^{1,3}$, 
V.V. Smirnova$^1$, V.A. Taperechkina$^2$\\
        (1) Institute for high energy physics, 142281 Protvino, Russia\\
        (2) Moscow State Academy of Instrument Engineering and \\
            Computer Science, Moscow, Russia\\
        (3) Serguei.Bitioukov@cern.ch, http://home.cern.ch/bityukov\\
       }

\maketitle

\begin{abstract}
\noindent
The reconstruction of the parameter of the model 
by the measurement of the random variable 
depending on this parameter is one of the main tasks of  
statistics. In the paper the notion of the statistically dual distributions 
is introduced. The approach, based on the 
properties of the statistically dual distributions,
to resolving of the given task is proposed. 
  

\end{abstract}
\thispagestyle{empty}

\section{Introduction}

As shown in refs.~(Bityukov~2000, Bityukov~2002, Bityukov~2003), 
in frame of frequentist approach we can construct the probability distribution 
of the possible magnitudes of the Poisson distribution parameter
to give the observed number of events $\hat n$ in Poisson stream of events.  
This distribution is
described as Gamma-distribution $\Gamma_{1, \hat n + 1}$ with
the probability density, which looks like Poisson distribution of 
probabilities, namely,

\begin{equation}
g_{\hat n}(\mu) = \displaystyle \frac{\mu^{\hat n}}{\hat n!} e^{-\mu},~ 
\mu > 0,~\hat n > -1,  
\end{equation}

\noindent
where $\mu$ is a variable and $\hat n$ is a parameter 
(in case of Poisson distributions $\mu$ is a parameter 
and $\hat n$ is a variable). It means, as shown below, that
we can estimate the value and error of Poisson distribution parameter by the 
measurement of mean value of the random variable of this Poisson distribution 
and by using the correspondent Gamma-distribution.

Let us name such distributions, which allow to exchange the 
parameter and the variable, conserving the same formula for the 
distribution of probabilities, as {\bf statistically dual distributions}.

Note, it is pure probabilistic (and, in this sense, frequentist) definition. 
As shown below, the properties of considered in this paper of two pairs 
of statistically dual distributions coincide with properties of conjugate 
families (which defined in frame of 
Bayesian Theory~(for example, Bernardo, 1994)).

In the next Section we show that Poisson and Gamma- distributions 
are statistically dual distributions and that normal distribution with 
constant variance is statistically self-dual distribution. 
Despite little differences in formulae, we also point 
on the interdependency between 
the negative binomial and $\beta-$ distributions, which allows to consider
these distributions as statistically quasi-dual distributions. 
The interrelation of the statistical duality and the estimation 
of the distribution parameters are discussed in Section~3.

\section{Statistically dual distributions}

Let $t(x,y)$ be a function of two variables. If the same function can 
be represented both as the distribution of the random variable $x$ with 
parameter $y$ and as the distribution of the random variable $y$ with 
parameter $x$, then these distributions can be named as 
{\bf statistically dual distributions}.

The {\bf statistical duality} of Poisson and Gamma-distributions follows from 
simple discourse.

Let us consider the Gamma-distribution with probability density

\begin{equation}
g_x(\beta,\alpha) = \displaystyle 
\frac{x^{\alpha-1}e^{-\frac{x}{\beta}}}{\beta^{\alpha}\Gamma(\alpha)}.   
\end{equation}

At change of standard designations of Gamma-distribution 
$\displaystyle \frac{1}{\beta}$, $\alpha$ and $x$ for  
$a$, $n+1$ and $\mu$ we get the following formula for probability density 
of Gamma-distribution 

\begin{equation}
g_n(a,\mu) = \displaystyle 
\frac{a^{n+1}}{\Gamma(n+1)} e^{-a\mu} \mu^{n},   
\end{equation}

\noindent
where $a$ is a scale parameter and $n + 1 > 0$ is a shape parameter. 
Suppose $a = 1$, then the formula of the probability density of 
Gamma-distribution $\Gamma_{1,n+1}$ is (here we repeat 
Eq.1)

\begin{center}
$g_n(\mu) = \displaystyle \frac{\mu^n}{n!} e^{-\mu},~ 
\mu > 0,~n > -1$.  
\end{center}

It is usually supposed that the probability of observing $n$ events in the 
experiment 
is described by Poisson distribution with parameter $\mu$, i.e.

\begin{equation}
f(n; \mu)  = \frac{{\mu}^n}{n!} e^{-\mu},~\mu > 0,~n \ge 0.
\end{equation}

One can see that the parameter and variable in Eq.1 and Eq.4 are exchanged
in other respects the formulae are identical.
As a result these distributions (Gamma and Poisson) are {\bf statistically 
dual distributions}. These distributions (Bityukov 2002)
are connected by identity~(see, also, identities in 
refs.~(Jaynes 1976, Frodesen 1979, Cousins 1995))

\begin{equation}
\displaystyle
\int_{\mu_1}^{\mu_2}{g_{m}(\mu) d\mu} + 
\sum_{i = n+1}^{m}{f(i;\mu_2)}~+
\displaystyle \int_{\mu_2}^{\mu_1}{g_{n}(\mu) d\mu} 
- \sum_{i = n+1}^{m}{f(i;\mu_1)} = 0~,
\end{equation}

\noindent
for any real $\mu_1 \ge 0$, $\mu_2 \ge 0$, and
integer $m > n \ge 0$.

The other example of statistically dual distribution is 
the normal distribution with mean value $a$ and constant variance $\sigma^2$:

\begin{equation}
\displaystyle \varphi(x;a,\sigma) = \displaystyle \frac{1}{\sqrt{2\pi}\sigma}
\displaystyle e^{-\frac{(x-a)^2}{2\sigma^2}},
\end{equation}

\noindent
where $x$ is real variable, $a$ and $\sigma > 0$ are real parameters.
Here we also can exchange the parameter $a$ and variable $x$ saving 
the formulized description of the probability density. It allows to estimate
the value of parameter $a$ by measurement of mean value for $x$. In this case
we must consider new density of probability with the real  
variable $a$ and real parameters $x$ and $\sigma > 0$

\begin{equation}
\displaystyle \phi(a;x,\sigma) = \displaystyle \frac{1}{\sqrt{2\pi}\sigma}
\displaystyle e^{-\frac{(x-a)^2}{2\sigma^2}}.
\end{equation}

\noindent
Such a way the normal distribution can be named as 
{\bf statistical self-dual distribution}.
The analogous identity (as Eq.5) for given distributions takes place

$$\displaystyle \int_{c}^{d}{\phi(a;c,\sigma)da} + 
\int_{c}^{d}{\varphi(x;d,\sigma)dx} + 
\int_{d}^{c}{\phi(a;d,\sigma)da} +
\int_{d}^{c}{\varphi(x;c,\sigma)dx} = 0$$

\noindent
or, simpler, 

\begin{equation}
\displaystyle \int_{c}^{d}{\phi(a;b,\sigma)da} - 
\int_{c}^{d}{\varphi(x;b,\sigma)dx}  = 0
\end{equation}

\noindent
for any real $b$, $c$ and $d$.


Let us present~(Eadie 1971) the probability density of $\beta-$distribution as 

\begin{equation}
\beta(x;n,m) = \frac{(n+m+1)!}{n!m!}x^n(1-x)^m,
\end{equation}

\noindent
with real variable $x$ ($0 < x < 1$) and integer non negative
parameters $n$ and $m$, and the probability distribution of negative 
binomial distribution as

\begin{equation}
P(k;n,p) = \frac{(n+k)!}{n!k!}p^{n+1}(1-p)^k,
\end{equation} 

\noindent
where $k$ is integer non negative variable, integer $n$ and real $p$ 
are parameters ($0 \le p \le 1$). 

The formulae, which  describe these distributions, are different, nevertheless,
as is mentioned in ref.~(Jaynes 1976), the negative binomial (Pascal) and 
$\beta-$ distributions have interrelation in the form of identity 

\begin{equation}
\int_0^p{\beta(x;n,m)dx} - \sum_{k=0}^m{P(k;n,p)} = 0.
\end{equation}

\noindent
Correspondingly, these distributions can be named as 
{\bf statistically quasi-dual distributions}.

\section{Statistical duality and estimation of the parameter of distribution} 

The identity (Eq.5) can be written in form~(Bityukov 2000, Bityukov 2003)

\begin{equation}
\displaystyle
\sum_{i = \hat n + 1}^{\infty}{f(i; \mu_1)} +
\int_{\mu_1}^{\mu_2}{g_{\hat n}(\mu) d\mu} + 
\sum_{i = 0}^{\hat n}{f(i;\mu_2)} = 1,~
\end{equation}

\noindent
i.e.

\begin{center}
$\displaystyle
\sum_{i = \hat n + 1}^{\infty}{\frac{\mu_1^ie^{-\mu_1}}{i!}} +
\int_{\mu_1}^{\mu_2}
{\frac{\mu^{\hat n}e^{-\mu}}{\hat n!}d\mu}
+ \sum_{i = 0}^{\hat n}{\frac{\mu_2^ie^{-\mu_2}}{i!}} = 1~$
\end{center} 

\noindent
for any real $\mu_1 \ge 0$ and $\mu_2 \ge 0$ and non negative integer 
$\hat n$. 

The definition of the confidence interval $(\mu_1,\mu_2)$
for Poisson distribution parameter $\mu$ as~(Bityukov 2000)

\begin{equation}
P(\mu_1 \le \mu \le \mu_2) = P(i \le \hat n|\mu_1) - P(i \le \hat n|\mu_2),
\end{equation}

\noindent 
where $\displaystyle 
P(i \le \hat n|\mu) = \sum_{i=0}^{\hat n}{\frac{\mu^ie^{-\mu}}{i!}}$,
allows to show
that a Gamma-distribution $\Gamma_{1,1+\hat n}$ is the probability distribution
of different values of $\mu$ parameter of Poisson distribution on condition 
that the observed value of the number of events is equal to $\hat n$. 
This definition is consistent with identity Eq.12 in contrast with another 
frequentist definitions of confidence interval (for example, if we suppose 
in Eq.12 that $\mu_1=\mu_2$ we have conservation of probability).
The right part of Eq.13 determines the frequentist sense of this definition.

Let us suppose that $g_{\hat n}(\mu)$ is the probability density of parameter 
of the Poisson distribution if number of observed events is equal to $\hat n$.
It is a conditional probability density. As it follows from formulae 
(Eqs.1,12), 
the $g_{\hat n}(\mu)$ is the density of Gamma-distribution by definition. 
Note, that the definition of conjugate families~
\footnote{Given a family $\cal F$ of pdf's $f(x|\theta)$ indexed by a 
parameter $\theta$, then a family, $\Pi$ of prior distributions is said to be 
{\it conjugate} for the family $\cal F$ if the posterior distribution 
of $\theta$ is in the family $\Pi$ for all $f \in \cal F$, all priors
$\pi(\theta)\in \Pi$ and all possible data sets $x$}
is based on using of 
probability density of parameter 
distribution~(Casella, 2001).

On the other hand: 
if $g_{\hat n}(\mu)$ is not equal to this probability density 
and the probability density of the Poisson parameter is the other 
function $h(\mu;\hat n)$ then there takes place another identity

\begin{equation} 
\displaystyle \sum_{i=\hat n+1}^{\infty} f(i;\mu_1) +
    \int_{\mu_1}^{\mu_2}{h(\mu;\hat n)d\mu} + 
\sum_{i=0}^{\hat n} f(i;\mu_2) = 1. 
\end{equation}

This identity is correct for any real $\mu_1 \ge 0$ and $\mu_2~\ge~0$ too.
The sums in the left part of this equation determine the boundary conditions 
on the confidence interval.

\bigskip

If we subtract Eq.14 from Eq.12 then we have

\begin{equation}
 \int_{\mu_1}^{\mu_2}{(g_{\hat n}(\mu) - h(\mu;\hat n))d\mu}  = 0.
\end{equation}

We can choose the $\mu_1$ and $\mu_2$ by the arbitrary way. Let us make this 
choice so that $g_{\hat n}(\mu)$ is not equal $h(\mu;\hat n)$ in the interval
$(\mu_1,\mu_2)$ and, for example, $g_{\hat n}(\mu) > h(\mu;\hat n)$ and
$\mu_2 > \mu_1$. In this case we have 

\begin{equation}
\int_{\mu_1}^{\mu_2}{(g_{\hat n}(\mu) - h(\mu;\hat n))d\mu} > 0
\end{equation}

\noindent
and we have contradiction (i.e. $g_{\hat n}(\mu)=h(\mu;\hat n)$ 
everywhere except, may be, a finite set of points).  
As a result we can mix Bayesian 
($g_{\hat n}(\mu)$) and frequentist ($f(k;\mu)$) probabilities 
without any logical inconsistencies. The identity (Eq.12) leaves no 
place for any prior except uniform ($\pi(\mu)=const$~\footnote{Bayesian 
methods are supposed that $\displaystyle g_{\hat n}(\mu)=
\frac{f(\hat n;\mu)\cdot\pi(\mu)}{\int{f(\hat n;\mu)\cdot\pi(\mu)d\mu}}$, 
where $\pi(\mu)$ is the prior probability density for $\mu$.}). 
As a result, we can construct the distribution of the errors
in the estimation of the parameter of Poisson distribution
by single measurement of a casual variable and, correspondingly, 
the confidence 
intervals, estimate the parameter by several measurements,
take into account systematics and statistical uncertainties of
measurements at statistical conclusions about the quality of planned 
experiments (Bityukov 2000, Bityukov 2003, Bityukov 2004). 

For the normal distribution the identity (Eq.8) can be written
in case of the observed value of random variable $\hat x$ as

\begin{equation}
\displaystyle \int_{-\infty}^{\hat x-c}{\varphi(x;\hat x,\sigma)dx} + 
\int_{\hat x-c}^{\hat x+d}{\phi(a;\hat x,\sigma)da} +
\displaystyle \int_{\hat x+d}^{\infty}{\varphi(x;\hat x,\sigma)dx} = 1
\end{equation}

\noindent
for any real $c$ and $d$ or, in analogy with Eq.12 for the case of 
Poisson-Gamma distributions,

\begin{equation}
\displaystyle \int_{\hat x}^{\infty}{\varphi(x;\hat x - c,\sigma)dx} + 
\int_{\hat x-c}^{\hat x+d}{\phi(a;\hat x,\sigma)da} +
\displaystyle \int_{-\infty}^{\hat x}{\varphi(x;\hat x+d,\sigma)dx} = 1
\end{equation}

\noindent
for any real $c \ge 0$ and $d \ge 0$.

These identities (18,19) also allow to say that conditional distribution 
(if observed value is $\hat x$) of true value of the parameter $a$
obeys normal distribution with mean value  $\hat x$ and 
constant variance $\sigma^2$
(i.e. here, in contrast with previous example, 
$\hat x$ is the unbiased estimator of the parameter $a$). 
As a result we can construct the distribution of the error and the confidence 
intervals of parameter $a$, take into account systematics and statistical 
uncertainties and so on 
in accordance with standard analysis of errors~(Eadie 1971).

So, the statistical duality allows to connect the estimation
of the parameter with the measurement of the random variable of the
distribution under study.

\section{Conclusion}

In the paper the notion of the statistically dual distributions is introduced.
The relation between the measurement of casual variable and the estimation
of the given distribution parameter for 
two pairs of  statistically dual distributions is presented. 
The proposed
approach allows to construct the distribution of measurement error of the
estimator of distribution parameter by the using of statistically dual
distribution.

Both considered cases of statistically
dual distributions 
(which are introduced in frame of frequentist approach), namely 
Poisson distribution versus Gamma distribution and Normal distribution 
versus Normal distribution, also belong to conjugate families
(which are defined in frame of Bayesian approach). 
For example~(Bityukov 2003 ACAT), the distributions conjugated to Poisson 
distributions were built by Monte Carlo method (i.e. in frame of frequentist
approach). The hypotheses testing confirms that these distributions are 
Gamma-distributions as expected in this case.

\subsubsection*{Acknowledgements}

The authors are grateful to V.A.~Kachanov, Louis
Lyons, N.V.~Krasnikov and 
V.F.~Obraztsov for the interest and useful comments, 
S.S.~Bityukov, 
R.D.~Cousins, Yu.P.~Gouz, G.~Kahrimanis, M.N.~Ukhanov 
for fruitful discussions and E.A.~Medvedeva for help in preparing
the paper. This work has been particularly supported 
by grants RFBR 04-01-97227 and RFBR 04-02-16020. 
 
\subsubsection*{References} 
 
J.M. Bernardo, A.F.M. Smith, {\it Bayesian Theory}, Wiley, Chichester, 1994. 

S.I.~Bityukov, N.V.~Krasnikov, V.A.~Taperechkina (2000).
Confidence intervals for Poisson distribution parameter.
{\it Preprint IFVE} 2000-61, Protvino;  {\it e-Print}: hep-ex/0108020, 2001. 

S.I.~Bityukov (2002). On the Signal Significance in the 
Presence of Systematic and Statistical Uncertainties. 
{\it Journal of high energy physics} {\bf 09}(060). 

S.I.~Bityukov and N.V.~Krasnikov (2003). Signal Significance in the Presence 
of Systematic and Statistical Uncertainties. {\it  Nuclear Instruments and 
Methods in Nuclear Research} {\bf A502}:795-798.  

S.I.~Bityukov, V.A.~Medvedev, V.V.~Smirnova, Yu.V.~Zernii (2003 ACAT). 
Experimental test of the probability density function 
of true value of Poisson distribution parameter by single observation of 
number of events. {\it  Nuclear Instruments and Methods in Nuclear Research} 
{\bf A534}:228-231; also e-Print: physics/0403069, 2004. 

S.I.~Bityukov and N.V.~Krasnikov (2004).  
The probability of making a correct decision in hypotheses testing 
as estimator of quality of planned experiments.  
in  G.~Erickson and Y.~Zhai (eds.), {\it
Bayesian Inference and Maximum Entropy Methods in Science and 
Engineering}, 23-th International Workshop on Bayesian Inference and 
Maximum Entropy Methods in Science and Engineering, Jackson Hole, 
Wyoming, 3-8 August 2003,  AIP Conference Proceedings, vol.{\bf 707}:455-464, 
Melville, NY, 2004; also e-Print:   physics/0309031 

See as an example,
G. Casella and R. L. Berger (2001),
{\it Statistical Inference} (2nd Edition), Duxbury Press.

R.D. Cousins (1995).  Why isn't every physicist a Bayesian ?
{\it Am.J.Phys} {\bf 63}:398-410.

See as an example, 
W.T.~Eadie, D.~Drijard, F.E.~James, M.~Roos, and B.~Sadoulet (1971).
{\it Statistical Methods in Experimental Physics,} North Holland,
Amsterdam.

A.G.~Frodesen, O.~Skjeggestad, H.~T$\o$ft (1979). 
{\it Probability and Statistics in Particle Physics.}
UNIVERSITETSFORLAGET, Bergen-Oslo-Troms$\o$, p.97.

E.T.~Jaynes (1976). Confidence intervals vs Bayesian Intervals.
In R.D.~Rosenkrantz (ed.), {\it E.T.~Jaynes: Papers on probability, 
statistics and statistical physics}, D.Reidel Publishing Company,
Dordrecht, Holland, 1983, p.165.

\end{document}